\documentclass[12pt]{article}
\usepackage[utf8]{inputenc}
\usepackage{enumerate}
\DeclareFontFamily{U}{mathb}{\hyphenchar\font45}
\DeclareFontShape{U}{mathb}{m}{n}{
<5> <6> <7> <8> <9> <10> gen * mathb
<10.95> mathb10 <12> <14.4> <17.28> <20.74> <24.88> mathb12
}{}

\DeclareSymbolFont{mathb}{U}{mathb}{m}{n}
\DeclareMathSymbol{\llcurly}{3}{mathb}{"CE}
\DeclareMathSymbol{\ggcurly}{3}{mathb}{"CF}
\DeclareFontFamily{U}{matha}{\hyphenchar\font45}
\DeclareFontShape{U}{matha}{m}{n}{
<5> <6> <7> <8> <9> <10> gen * matha
<10.95> matha10 <12> <14.4> <17.28> <20.74> <24.88> matha12
}{}
\DeclareSymbolFont{matha}{U}{matha}{m}{n}
\newenvironment{E}{\begin{equation}}{\end{equation}}
\def\proof{\noindent{\bf Proof: }}
\def\qed{ \hskip 20pt{\vrule height7pt width6pt depth0pt}\hfil}
\def\bs{\llcurly}
\def\forb{{\mathrm{forb}}}
\def\bsf{{\mathrm{ Bh}}}
\def\ex{{\mathrm{ex}}}
\def\Av{{\mathrm{Avoid}}}
\def\avbh{{\mathrm{BAvoid}}}
\def\bAv{{\mathrm{BAvoid}}}

\def\0{{\bf 0}}
\def\1{{\bf 1}}

\newcommand{\ncols}[1]{\| #1 \|}
\newcommand{\rf}[1]{(\ref{#1})}
\newcommand{\drf}[1]{Definition~\ref{#1}}
\newcommand{\trf}[1]{Theorem~\ref{#1}}
\newcommand{\lrf}[1]{Lemma~\ref{#1}}

\newcommand{\corf}[1]{Conjecture~\ref{#1}}

\newcommand{\srf}[1]{Section~\ref{#1}}

\newtheorem{thm}{Theorem}[section]
\newtheorem{lemma}[thm]{Lemma}

\newtheorem{cor}[thm]{Corollary}
\newtheorem{conj}[thm]{Conjecture}
\newtheorem{defn}[thm]{Definition}
\newtheorem{remark}[thm]{Remark}

\renewcommand{\l}{\ell}
\usepackage{amssymb}
\usepackage{amsmath}
\title{ Forbidden Berge hypergraphs}
\author{R.P. Anstee\thanks{Research supported in part by NSERC}, Santiago Salazar\thanks{Research supported in part by NSERC USRA} \\Mathematics Department\\The University of British Columbia\\Vancouver, B.C. Canada V6T 1Z2\\ {\small\texttt{anstee@math.ubc.ca}},  {\small\texttt{santiago.salazar@zoho.com}}
\\\mbox{\ }}

\begin{document}
\maketitle
\begin{abstract}
A \emph{simple} matrix is a (0,1)-matrix with no repeated columns. For a (0,1)-matrix $F$, we say
 that a (0,1)-matrix $A$ has $F$ as a \emph{Berge hypergraph} if there is a submatrix $B$ of $A$ and some row and column permutation of 
 $F$, say $G$, with $G\le B$.   Letting $\ncols{A}$ denote the number of columns in $A$, we define the extremal function $\bsf(m,{ F})=\max\{\ncols{A}\,:\, A \hbox{ is }m\hbox{-rowed simple matrix with no Berge hypergraph }F\}$. We determine the asymptotics of $\bsf(m,F)$ for all $3$- and $4$-rowed $F$ and most $5$-rowed $F$.  For certain $F$, this becomes the problem of determining the maximum number of copies of $K_r$ in a $m$-vertex graph that has no $K_{s,t}$ subgraph, a problem studied by Alon and Shinkleman.

\end{abstract}

\section{Introduction}

This paper  explores  forbidden Berge hypergraphs and their relation to forbidden configurations. 
Define a matrix to be \emph{simple} if it is a (0,1)-matrix with no repeated columns. Such a matrix can be viewed as an element-set incidence matrix. 
Given two (0,1)-matrices $F$ and $A$, we say $A$ has $F$ as a \emph{Berge hypergraph} and write $F\bs A$  if there is a submatrix $B$ of $A$ and a row and column permutation of $F$, say $G$, with $G\le B$. The paper of Gerbner and Palmer \cite{GP} introduces this concept to generalize the notions of Berge cycles and Berge paths in hypergraphs. Let $F$ be $k\times \ell$. A Berge hypergraph associated with the object $F$ is a hypergraph whose restriction to a set of $k$ elements yields a hypergraph that `covers' $F$.
Berge hypergraphs are related to the notion of a \emph{pattern} $P$ in a (0,1)-matrix $A$ which has been extensively studied and is quite challenging \cite{FH}.  
We say $A$ has pattern $P$ if there is a submatrix $B$ of $A$ with $P\le B$. The award winning result of Marcus and Tardos \cite{MT}  concerns avoiding a pattern corresponding to a permutation matrix. Row and column order matter to patterns.

We use heavily the concept of  a \emph{ configuration}; see \cite{survey}.  We say $A$ has a configuration $F$ if there is a submatrix $B$ of $A$ and a row and column permutation of $F$, say $G$, with $B= G$. Configurations care about the 0's as well as the 1's in $F$ but do not care about row and column order.  In set terminology the notation \emph{trace} can be used. 

 For a subset of rows $S$, define $A|_S$ as the submatrix of $A$ consisting of rows $S$ of $A$.  Define $[n]=\{1,2,\ldots ,n\}$. If $F$ has $k$ rows and $A$ has $m$ rows and $F\bs A$ then there is a $k$-subset $S\subseteq[m]$ such that $F\bs A|_S$. For two $m$-rowed matrices $A,B$, use $[A\,|\,B]$ to denote the concatenation of $A,B$ yielding a larger $m$-rowed matrix. Define 
 $t\cdot A=[A\,A\,\cdots A]$ as the matrix obtained from concatenating $t$ copies of $A$. 
  Let $A^c$ denote the (0,1)-complement of $A$. Let $\1_a\0_b$ denote the $(a+b)\times 1$ vector of $a$ 1's on top of $b$ 0's. We use $\1_a$ instead of $\1_a\0_0$. Let $K_k^{\ell}$ denote the $k\times \binom{k}{\ell}$ simple matrix of all columns of $\ell$ 1's on $k$ rows and let $K_k=[K_k^0K_k^1K_k^2\cdots K_k^k]$.

Define $\ncols{A}$ as the number of columns  of $A$. Define our extremal problem as follows:
$$\avbh(m,{\cal F})=\{A\,:\,A\hbox{ is $m$-rowed, simple, }F\not\llcurly A\hbox{ for all }F\in{F}\},$$ 
$$\bsf(m,{\cal F})=\max_A\{\ncols{A}\,:\,A\in\avbh(m,{\cal F})\}.$$
We are mainly interested in ${\cal F}$ consisting of a single forbidden Berge hypergraph $F$.
When $|{\cal F}|=1$ and ${\cal F}=\{F\}$, we write $\avbh(m,F)$ and $\bsf(m,F)$.

The main goal of this paper is to explore the asymptotic growth rate of $\bsf(m,F)$ for a given $k\times \l$ $F$. \trf{classifyk=3} handles $k=3$, \trf{classifyk=4} handles $k=4$ and \trf{classifyk=5} handles $k=5$ (modulo \corf{conjC4}).  The results apply some of the proof techniques (and results)  for Forbidden configurations \cite{survey}.  We have some interesting connections with $\ex(m,K_{s,t})$ (the maximum number of edges in a graph on $m$ vertices with no complete bipartite graph $K_{s,t}$ as a subgraph) and $\ex(m,K_n,K_{s,t})$ \cite{AS} (the maximum number of subgraphs $K_n$ in a graph on $m$ vertices with no complete bipartite graph $K_{s,t}$ as a subgraph). Two such results are \trf{K2t}  and \trf{I3Ik}. We also obtain in \trf{tree}, that if $F$ is the vertex-edge incidence matrix of a tree $T$, then $\bsf(m,F)$ is $\Theta(m)$ analogous to $\ex(m,T)$.  Note that $K_k$ has two meanings in this paper that are hopefully clear by context namely as the complete graph on $k$ vertices or as  the matrix $[K_k^0K_k^1K_k^2\cdots K_k^k]$.

We first make some easy observations.

\begin{remark}Let $F,F'$ be two $k\times\ell$ (0,1)-matrices satisfying $F\bs F'$.  Then $\bsf(m,F)\le \bsf(m,F')$.\end{remark}

The related extremal problem for forbidden configurations is as follows:
$$\Av(m,{\cal F})=\{A\,:\,A\hbox{ is $m$-rowed, simple,}F\not\prec A\hbox{ for all }F\in{F}\},$$ 
$$\forb(m,{\cal F})=\max_A\{\ncols{A}\,:\,A\in\Av(m,{\cal F})\}.$$
When $|{\cal F}|=1$ and ${\cal F}=\{F\}$, we write $\Av(m,F)$ and $\forb(m,F)$. 
There are striking differences between $\bsf(m,F)$ and $\forb(m,F)$ such as $\trf{tree}$ for Berge hypergraphs and $\trf{treeforb}$ for Forbidden configurations.
Note that the two notions of Berge hypergraphs and Configurations coincide when $F$ has no 0's.

\begin{remark}\label{J} Let $F$ be a (0,1)-matrix. Then $ \bsf(m,F)\le\forb(m,{ F}) $. If $F$ is  a matrix of 1's  
then  $\bsf(m,F)=\forb(m,{ F})$.\end{remark}

Note that any forbidden Berge hypergraph $F$ can be given as a family ${\cal B}(F)$ of forbidden configurations by replacing the 0's of $F$ by 1's in all possible ways.  
Define
\begin{E}{\cal B}(F)=\{B \hbox{ is a (0,1)-matrix} \,:\,F\le B\}.\label{bergefamily}\end{E}
Isomorphism can reduce the required set of matrices to consider, for example ${\cal B}({I_2})$ which has 4 matrices  satisfies:
$$\avbh(m,{\cal B}({I_2})=\avbh(m,
\left\{
\left[\begin{array}{cc}1&0\\0&1\\\end{array}\right],
\left[\begin{array}{cc}1&1\\0&1\\\end{array}\right],
\left[\begin{array}{cc}1&1\\1&1\\\end{array}\right]
\right\}).$$
\begin{remark} $\bsf(m,F)=\forb(m,{\cal B}(F))$. \end{remark}

A product construction is helpful here. Let $A$, $B$  be  $m_1\times n_1$ and $m_2\times n_2$ matrices respectively. We define $A\times B$ as the $(m_1+m_2)\times n_1n_2$ matrix whose columns are obtained by placing a column of $A$ on top of a column of $B$ in all $n_1n_2$ possible ways. This extends readily to $p$-fold products. Let $I_t=K_t^1$ denote the $t\times t$ identity matrix.  In what follows you may assume $p$ divides $m$ since we are only concerned with asymptotic growth with respect to $m$
$$\hbox{The  }p\hbox{-fold product }
\overbrace{I_{m/p}\times I_{m/p}\times \cdots\times I_{m/p}}^p,$$
 is an $m\times m^p/p^p$ simple matrix.  This corresponds to the vertex-edge incidence matrix of the complete $p$-partite hypergraph with parts $V_1,V_2,\ldots,V_p$ each of size $m/p$ so that
 $\{v_1,v_2,\ldots ,v_p\}$ is an edge if and only if $v_i\in V_i$ for $i=1,2,\ldots ,p$.   These products sometimes yield the asymptotically best (in growth rate) constructions avoiding $F$ as a Berge hypergraph.

 \begin{remark} Let $F$ be a given $k\times \ell$ (0,1)-matrix so that $F\not\bs I_{m/p}\times I_{m/p}\times\cdots\times I_{m/p}$ (a $p$-fold product). Then  $\bsf(m,F)$ is $\Omega(m^{p})$.  \end{remark} 
  
  Sometimes the product may contain the best construction using  the following idea from  \cite{AKRS}, that when given two matrices $F,P$ where $P$ is $m$-rowed then 
$$f(F,P)=\max_A\{\ncols{A}\,|\, A\hbox{ is }m{-rowed}, A\bs P\hbox{ and } F\not\bs A\}.$$
Thus \trf{C4} yields $\bsf(m,I_2\times I_2)$ is $\Theta(f(I_2\times I_2,I_{m/2}\times I_{m/2}))$.
The result \lrf{K3inKst} indicates that things must be more complicated for general $s,t$.

 A shifting argument works nicely here. 
 We let $T_i(A)$ denote the matrix obtained from $A$ by attempting to replace 1's in row $i$ by 0's . We do not replace a 1 by a 0 in row $i$ and column $j$ if the resulting column is already present in $A$ otherwise we do replace the 1 by a 0.  We have that $\ncols{T_i(A)}=\ncols{A}$ and if $A$ is simple then  $T_i(A)$ is simple.  
 \begin{lemma}  Given $A\in\avbh(m,F)$, there exists a matrix $T(A)\in\avbh(m,F)$ with $\ncols{A}=\ncols{T(A)}$ and $T_i(T(A))=T(A)$ for $i=1,2,\ldots ,m$.  \label{shifting}\end{lemma}
 
 \proof It is automatic that $\ncols{A}=\ncols{T_i(A)}$. We note that $F\not\llcurly A$ implies $F\not\llcurly T_i(A)$.  
 Replace $A$ by $T_i(A)$ and repeat.  Let $T^*(A)=T_{m}(T_{m-1}(\cdots T_1(A)\cdots))$. Either $T^*(T^*(A))$ contains 
 fewer 1's than $T^*(a)$ or we have $T_i(T^*(A))=T^*(A)$ for $i=1,2,\ldots ,m$.  In the former case replace $A$ by $T^*(A)$ and repeat. In the latter case let $T(A)=T^*(A)$.  Since the number of 1's in $A$ is finite, then the algorithm will terminate with our desired matrix $T(A)$. 
\qed
 
 \vskip 10pt
 Typically $T(A)$ is referred to as a \emph{downset} since when the columns of $T(A)$ are interpreted as a set system ${\cal T}$ then if $B\in{\cal T}$ and $C\subset B$ then $C\in{\cal T}$.  Note that if $T(A)$ has a column of sum $k$ with 1's on rows $S$, then 
 $K_k\bs T(A)|_S$ and moreover the copy of $K_k$ on rows $S$ can be chosen with 0's on all other rows. An easy consequence is that for $A\in\avbh(m,F)$ where $F$ is $k$-rowed and simple then we may assume $A$ has no columns of sum $k$.

\vskip 10pt
\section{General results}

This section provides a number of results about Berge hypergraphs that are used in the paper. The following results from forbidden configurations were useful.

\begin{thm}\cite{AF} Let $k,t$ be given with $t\ge 2$. Then $\forb(m,t\cdot \1_k)=\forb(m,t\cdot K_k)$ and is $\Theta(m^{k})$.  \label{block1s}\end{thm}

\begin{thm}\cite{AS} Let $k,t$ be given. Then $\forb(m,[\1_k\,|\,t\cdot K_k^{k-1}])$ is $\Theta(m^{k-1})$. \label{boundary}\end{thm}

\begin{thm}\cite{AFl} Let $F$ be a $k$-rowed simple matrix. 
Assume there is some pair of rows $i,j$ so than no column of $F$ contains
0's on rows $i,j$, there is some pair of rows $i,j$ so than no column of $F$ contains
1's on rows $i,j$ and there is some pair of rows $i,j$ so than no column of $F$ contains
$I_2$ on rows $i,j$. Then $\forb(m,F)$ is $O(m^{k-2})$. \label{smallboundary}\end{thm}

\begin{defn} Let $F$ be a $k$-rowed (0,1)-matrix. Define $G(F)$ as the graph on 
$k$ vertices such that we join vertices $i$ and $j$ by an edge if and only if there is a column in $F$ with 1's in rows $i$ and $j$. 
Let $\omega(G(F))$ denote the size of the largest clique in $G(F)$ and $\chi(G(F))$ is the chromatic number of $G(F)$. Let $\alpha(G(F))$ denote the size of the largest independent set
in $G(F)$. \label{G(F)}\end{defn}

 \begin{lemma} Let $F$ be given.  Then $\bsf(m,F)$ is $\Omega(m^{\chi(G(F))-1})$ and hence $\Omega(m^{\omega(G(F))-1})$.
 \label{X(F)}\label{clique}\end{lemma}
 
 \proof Let $p=m/(\chi(G(F))-1)$.  Let 
 $$A=\overbrace{ I_p\times I_p\times\cdots\times  I_p}^{\chi(G(F))-1}.$$
 If $F\llcurly A$ and the rows of the $\chi(G(F))-1$ fold product containing  $F$ are $T$ then we obtain $\chi(G(F))-1$ disjoint sets $T_1,T_2,\ldots ,T_{X(F)-1}$ with 
 $T_i=T\cap \{(i-1)p+1,(i-1)p+2,\ldots ,ip\}$ and $A|_{T_i}\bs I_{|T_i|}$.  This contradicts the  definition of $\chi(G(F))$ 
 and so $F\not\llcurly A$. Thus  $\bsf(m,F)$ is $\Omega(m^{\chi(G(F))-1})$.
  Note that $\chi(G)\ge \omega(G)$. 
 \qed

\begin{lemma}\label{construction} If  $2\cdot\1_t\bs F$ then  $\bsf(m,F)$ is $\Omega(m^t)$ \end{lemma}
\proof  $F$ is not  a Berge hypergraph of the $t$-fold product $I_{m/t}\times I_{m/t}\times\cdots\times I_{m/t}$. \qed

\begin{thm} Let $k$ be given and assume $m\ge k-1$. Then $\bsf(m,I_k)=2^{k-1}$.
\label{Ik}\end{thm}
\proof  The construction  consisting of $K_{k-1}$ with $m-k+1$ rows of 0's added yields $\bsf(m,I_k)\ge 2^{k-1}$. 
The largest $m$-rowed matrix which avoids $I_1=[1]$ as a Berge hypergraph is $[\0_m]$ since it avoids $I_1$ as Berge hypergraphs. This proves the base case $k=1$ and the following is the inductive step. 

Let $A\in\avbh(m,I_k)$.  Let $B$ be obtained from $A$ by removing any rows of 0's  so that $B$ is simple and every row of $B$ contains a 1. If $B$ has $k-1$ rows then $\ncols{A}=\ncols{B}\le 2^{k-1}$ which is our bound. Assume $B$ has at least $k$ rows.  Either $\ncols{B}\le 2^{k-1}$ in which case we are done or $\ncols{B}>2^{k-1}>2^{k-2}$ and so  by induction, $B$ must contain $I_{k-1}$ as a Berge hypergraph.  Permute $B$ to form the block matrix 
$$B=\left[\begin{array}{@{}c@{}|@{}c@{}}\begin{array}{c}C\\ \hline E\\ \end{array}&\begin{array}{c}D\\ \hline G\\ \end{array}\\ \end{array}\right]$$
where $C$ is $(k-1)\times (k-1)$ with $I_{k-1}\bs C$.  Then $G$  must be the matrix of 0's or else $I_{k}\bs B$. Thus $D$ is simple.   Since all rows of $B$  contain a 1,  then $E$ must have a 1.  If $E$ contains a 1 then $I_{k-1}\not\llcurly D$ and so $\ncols{D}\le 2^{k-2}$. This gives  $\ncols{B}=\ncols{C}+\ncols{D}= k-1 + 2^{k-2}\le 2^{k-1}$. Thus
 $\ncols{A}=\ncols{B}\le 2^{k-1}$. \qed
\vskip 10pt

While \trf{Ik} establishes a constant bound for the Berge hypergraph $I_k$, we can see that this follows from a result of Balogh and Bollob\'as  \cite{BB}. Let $I_k^c=K_k^{k-1}$ denote the $k\times k$  (0,1)-complement of $I_k$ and let $T_k$ denote the $k\times k$ upper triangular  (0,1)-matrix with a 1 in row $i$ and column $j$ if and only if $i\le j$.

\begin{thm} \cite{BB} Let $k$ be given. Then there is a constant {$c_k$} so that
\break $\forb(m,\{I_k,I_k^c,T_k\})={c_k}$.\label{BB}\end{thm}
A corollary of Koch and the first author \cite{AKoch} gives one way to apply this result.
.
\begin{thm} \cite{AKoch} Let ${\cal F}=\{F_1,F_2,\ldots ,F_t\}$ be given. There are two possibilities. Either $\forb(m,{\cal F})$ is $\Omega(m)$ or there exist 
$\ell,i,j,k$  with $F_i\prec I_{\ell}$,  with $F_j\prec I_{\ell}^c$ and  with $F_k\prec T_{\ell}$ in which case  there is a constant $c$ with $\forb(m,{\cal F})= c$.\label{classify}\end{thm}

We apply this result to a forbidden Berge hypergraph $F$  using the family ${\cal B}(F)$ from \rf{bergefamily} which contains the $k\times \ell$ matrix of 1's.  Noting that $I_{k+\ell+1}^c$ contains a $k\times \ell$ block of 1's  and  $T_{k+\ell}$ contains a $k\times \ell$ block of 1's we obtain the following.  

\begin{cor} Let $F$ be a $k\times \ell$ (0,1)-matrix.  Then either $\bsf(m,F)$ is $\Omega(m)$ or $F\bs I_{k+\ell}$ in which case 
$\bsf(m,F)$ is $O(1)$. \label{constantlinear}\end{cor}

 The following Lemma (from standard induction in \cite{survey}) was quite useful for Forbidden Configurations.
 
 \begin{lemma} Let $F$ be a $k\times \ell$ (0,1)-matrix and let $F'$ be a $(k-1)\times \ell$ submatrix of $F$.
 Then $\bsf(m,F)=O(m\cdot\bsf(m,F'))$. \label{induction}\end{lemma} 
 \proof Let $A\in\avbh(m,F)$. If we delete row 1 of $A$, then the resulting matrix may have columns that appear twice. We may permute the columns of $A$ so that
 $$A=\left[\begin{array}{cc}0\,0\cdots 0&1\,1\cdots 1\\ B\quad C&C\quad D\end{array}\right],$$
 where $[B\,C\,D]$ and $C$ are simple $(m-1)$-rowed matrices. We have  $[B\,C\,D]\in\avbh(m-1,F)$ and $C\in\avbh(m-1,F')$ (if $F'\bs C$ then $F\bs A$). Then 
 $$\ncols{A}=\ncols{[BCD]}+\ncols{C}\le \bsf(m-1,F)+\bsf(m-1,F'),$$
 which yields the desired bound by induction on $m$. \qed

 \begin{lemma}
	\label{rowsumik}
	Let A be a $k$-rowed (0,1)-matrix, not necessarily simple, with all row sums at least  $kt$. Then $t\cdot I_k \bs A$.
\end{lemma}
\proof
We use induction on $k$ where the case $k=1$ and $I_1=[1]$ is easy.
Choose $t$ columns from $A$ containing a 1 in row 1 and remove them and row 1 resulting in a matrix $A'$. The row sums of $A'$ 
will be at least $(k-1)t$ and so we may apply induction. Thus $(t-1)\cdot I_k \bs A'$ and so we obtain $t\cdot I_k \bs A$. \qed
\vskip 10pt

An interesting corollary is that if we have an $m$-rowed matrix $A$ with all rows sums at least $kt$ then $t\cdot I_k\bs A|_S$  for all $S\in\binom{[m]}{k}$.

\begin{lemma} Let $A$ be a given  $m$-rowed matrix and let ${\cal S}$ be a family of subsets of $[m]$ with the property that $|S|\le k$ for  all $S\in{\cal S}$.  Let $c$ be given. Then by deleting at most $c\left(\binom{m}{k}+\binom{m}{k-1}+\cdots+\binom{m}{1}\right)$ columns from $A$ we can obtain a matrix $A'$ so that for each $S\in{\cal S}$, $A'|_S$ either has more than $c$ columns with  1's on all the rows of $S$ or has no columns with 1's on all the rows of $S$. \label{longsupply}\end{lemma}
\proof  For each subset of $S\in{\cal S}$, if the number of columns of $A|_{S}$ with 1's on the rows of $S$ is at most $c$, then delete all such columns. Repeat. The number of deleted columns is at most 
$\sum_{S\in{\cal S}}c\le c\left(\binom{m}{k}+\binom{m}{k-1}+\cdots+\binom{m}{1}\right)$. \qed
\vskip 10pt

\begin{lemma}(Reduction Lemma) Let $F=[G\,|\,t\cdot [H\,K]]$. Assume $H,K$ are simple and have column sums at most $k$. 
Also assume  for each column $\alpha$ of $K$, there is a column $\gamma$ of $[G\,H]$ with $\alpha\le\gamma$. 
Then there is a constant $c$ so that
$\bsf(m,F)\le\bsf(m,[G\,H])+cm^{k}$. \label{morereduction}\end{lemma}

\proof  We let $A\in\avbh(m,[G\,|\,t\cdot [H\,K]])$ and $c=\ncols{G}+t\ncols{H}+t\ncols{K}$. Applying \lrf{shifting}, assume $T_i(A)=A$ for all $i$ and so, when columns are viewed as sets, the columns form a downset.  Form ${\cal S}$ as the union of all  $S$ so that $[H\,K]$ has a column with 1's on the rows $S$. Then, applying \lrf{longsupply}, delete at most $cm^{k}$ columns to obtain a matrix $A'$. Now if $[G\,H]\bs A'$ on rows $S$, then each column contributing to $H$ will appear $c$ times in $A'|_S$.  

 Moreover  each column $\gamma$ of $G$ will appear at least $c$ times in $A'|_S$ and so if $\alpha$ is a column of $K$ and  $\gamma$ is a column of $G$ with $\alpha\le\gamma$, then we have $t\cdot\alpha\bs t\cdot\gamma$. Hence 
$[G\,|\,t\cdot [H\,K]]\bs A|_S$, a contradiction.   The choice of $c$ is  required, for example, when the columns contributing to $[G\,H]$ all have $A|_S=\1$. \qed

\vskip 10pt
This is reminiscent of \lrf{reduction} but the rules for eliminating columns of small column sum (at most $k$) are slightly more strict.
The following are two important applications.
We use the notation $K_p\backslash \1_p$ to denote the matrix obtained from $K_p$ by deleting the column of $p$ 1's. 

\begin{thm}\label{1pxIk-p}
Let $H(p,k,t)=[1_p\times I_{k-p}\,t\cdot[\,\1_p\times\0_{k-p}\,|\,(K_p\backslash \1_p)\times [\0_{k-p}\,\,I_{k-p}]]$, namely
\begin{E}H(p,k,t)=  \left[\begin{array}{cccc}1&1&\cdots&1\\ \vdots&\vdots&\cdots&\vdots\\ 1&1&\cdots&1\\ 1&0&0&0\\ 0&1&0&0\\ 0&0&\ddots&0\\ 0&0&0&1\\ \end{array}\,\,\,
t\cdot\left[
\begin{array}{c}1\\ \vdots\\ 1\\ 0\\ 0\\ \vdots\\ 0\\ \end{array}
\,
\begin{array}{c}  \\ K_p\backslash\1_p\\   \times \\
\left[\begin{array}{ccccc}0&1&0&0&0\\ 0&0&1&0&0\\ 0&0&0&\ddots&0\\ 0&0&0&0&1\\ \end{array}\right]\\
\end{array}
\right]\right] .\label{H(p,k,t)}\end{E}
Then $\bsf(m,H(p,k,t))$ is $\Theta(m^p)$. Moreover if we add to $H(p,k,t)$ any column not already present $t$ times in $H(p,k,t)$ to obtain  $F'$,
 then $\bsf(m,F')$ is $\Omega(m^{p+1})$.
\end{thm}
\proof take $F=H(p,k,t)$. Given that $F$ has a column of $p+1$ 1's then $\omega(G(F))\ge p+1$ and so \lrf{clique} yields $\bsf(m,F)$ is $\Omega(m^p)$.

 To apply  Reduction \lrf{morereduction}, set $F=[G\,| \,t\cdot [H\,K]]$ with $G$ to be the first $k-p$ columns of $F$ and with $K$ to be the remaining $1+(2^p-1)\times (k-p)$ columns of $F$ when $t=1$ and with $H$ absent.  Now $\bsf(m,F)\le\bsf(m,G)+cm^{p}$ for $c=\ncols{G}+t\ncols{K}$. Applying \lrf{induction} repeatedly (in essence deleting the first $p$ rows of $G$) we obtain $\bsf(m,G)=O(m^k\bsf(m,I_{k-p}))$ and so with
 \lrf{Ik} this yields $\bsf(m,G)$ is $O(m^p)$. Then  $\bsf(m,H(p,k,t))$ is $\Theta(m^p)$. 
 
 The remaining remarks concerning adding a column to $H(p,k,t)$  are covered in \lrf{generalmaximal}.
\qed

\vskip 10pt
Note that $\bsf(m,H(k-1,k,t))$ follows from \trf{boundary}.  There is a more general form of $H(p,k,t)$ as follows. 

\begin{defn}
Let $A$ be a given (0,1)-matrix. Let ${\cal S}(A)$ denote the matrix of 
all columns $\alpha$ so that there exists a column $\gamma$ of $A$ with $\alpha\le\gamma$ and $\alpha\ne\gamma$. \qed \end{defn}
Let 
\begin{E}H((a_1,a_2,\ldots,a_s),t)=[I_{a_1}\times I_{a_2}\times\cdots\times I_{a_s}\,|\,t\cdot{\cal S}([I_{a_1}\times I_{a_2}\times\cdots\times I_{a_s}])]\label{generalH}\end{E}
Then $H(p,k,t)$ is $H((a_1,a_2,\ldots,a_s),t)$ where $s=p+1$ and $a_1=a_2=\cdots=a_p=1$ and $a_{p+1}=k-p$. The upper bounds 
of \trf{1pxIk-p} do not generalize but the second part of the proof continues to hold.

\begin{lemma}Let $H((a_1,a_2,\ldots,a_s),t)$ be defined as in \rf{generalH}. Then \hfil\break $\bsf(m,H((a_1,a_2,\ldots,a_s),t))$ is 
$\Omega(m^{s-1})$.  Moreover if we add to $H((a_1,a_2,\ldots,a_s),t)$ any column $\alpha$ not already present $t$ 
times in $H((a_1,a_2,\ldots,a_s),t)$ then 
$\bsf(m,[H((a_1,a_2,\ldots,a_s),t)\,|\,\alpha])$
 is $\Omega(m^{s})$. \label{generalmaximal}\end{lemma}
 \proof  The lower bound for $\bsf(m,H((a_1,a_2,\ldots,a_s),t))$ follows from $(s-1)$-fold product
 $I_{m/(s-1)}\times I_{m/(s-1)}\times\cdots\times I_{m/(s-1)}$ since $H((a_1,a_2,\ldots,a_s),t))$ has columns of sum $s$.
 
 There are two choices for $\alpha$. First we can choose $\alpha$ to be a column 
 in $I_{a_1}\times I_{a_2}\times\cdots\times I_{a_s}$ and so $\alpha$ has $s$ 1's.  Then $2\cdot\1_{s}\bs [\alpha\alpha]$ so that $\bsf(m,[\alpha\alpha]$ is $\Theta(m^{s})$ by \trf{block1s}. 
 
Second choose $\alpha$ to be a column not already present in $H((a_1,a_2,\ldots,a_s),t)$.  Let $G=G(H((a_1,a_2,\ldots,a_s),t))$ be the graph defined in \drf{G(F)} on $a_1+a_2+\cdots +s_s$ vertices corresponding to rows of $H((a_1,a_2,\ldots,a_s),t)$.   Our choice of 
 $\alpha$ has a pair of rows $h,\ell$ so that $\alpha$ has 1's in both rows $h$ and $\ell$ and the edge $h,\ell$ is not in $G$. We deduce that  
$[H((a_1,a_2,\ldots,a_s),t)\,|\,\alpha]$ has $s+1$ rows $S$ such that for every pair $i,j\in S$, there is a column with 1's in both rows $i$ and $j$, i.e. $G(H((a_1,a_2,\ldots,a_s),t))$ has a clique of size $s+1$. Thus by \lrf{clique},  
$\bsf(m,[H((a_1,a_2,\ldots,a_s),t)\,|\,\alpha])$
 is $\Omega(m^{s})$. \qed
 
 \vskip 10pt  Thus all but the upper bounds for \trf{1pxIk-p} follow from \lrf{generalmaximal}.
 The following application  requires \corf{conjC4} to be true. Note that $\1_1\times C_4$ is $I_1\times I_2\times I_2$.
\begin{thm}\label{1xC4}
Assume $\bsf(m,\1_1\times C_4)$ is $\Theta(m^2)$. Then  $\bsf(m,H((1,2,2),t))$ is $\Theta(m^2)$. Moreover if we add to $H((1,2,2),t)$ any column $\alpha$ not already present $t$ times in  \hfil\break$H((1,2,2),t)$ to obtain $[H((1,2,2),t)\,|\,\alpha]$,
 then $\bsf(m,[H((1,2,2),t)\,|\,\alpha])$ is $\Omega(m^3)$.
\end{thm}
\proof Take $G=\1_1\times I_2\times I_2=\1_1\times C_4$ and take $K$ to be the remainder of the columns of $H((1,2,2),1)$ and then  apply Reduction \lrf{morereduction} and the hypothesis that $\bsf(m,\1_1\times C_4)$ is $\Theta(m^2)$ to obtain the upper bound. 

The rest follows from \lrf{generalmaximal}. \qed

\vskip 10pt

The following monotonicty result seems obvious but note that monotonicity is only conjectured to be true for forbidden configurations. 
\begin{lemma}\label{monotonicity} Assume $F$ is a $k\times \ell$ matrix and assume $m\ge k$, Then
$\bsf(m,F)\ge\bsf(m-1,F)$. \end{lemma}
\proof Let $F'$ be the matrix obtained from $F$ by deleting rows of 0's, if any. Then for $m\ge k$,
$A\in\avbh(m,F)$ if and only if $A\in\avbh(m,F')$. Now assume $A\in\avbh(m,F')$ with $m\ge k$. Then form $A'$ from $A$ by adding a single row or 0's. Then $A'\in\avbh(m+1,F')$ 
with $\ncols{A}=\ncols{A'}$.\qed
\vskip 10pt

The following allows  $F$ to have rows of 0's which  contrasts with Reduction \lrf{morereduction}.

\begin{lemma}
	\label{reduction}
	Let $F$ be a $k\times \ell$ matrix. Then $\bsf(m,[F\,|\,t\cdot I_k]) \leq \bsf(m,F) + (tk+\ell)m$.
\end{lemma}
\proof
Let $A \in \bAv(m,[F\,|\,t\cdot I_k]$. For any row in $A$ of row sum $r$ we may remove that row and the $r$ columns containing a 1 on that row and the remaining $(m-1)$-rowed matrix is simple. In this way remove all rows with row sum at most $ tk+l$ and call the 
remaining simple matrix $B$ and assume it has $m'$ rows. Then $\ncols{A}\le\ncols{B}+(tk+\ell)(m-m')$. Suppose $B$ contains $F$ on some $k$-rows $S\subseteq \binom{[m']}{k}$. Remove the columns containing $F$ from $B$ to obtain $B'$ and now the  rows of $B'$  have row sum $\geq tk$. By \lrf{rowsumik}, $t \cdot I_k$ is contained in $B'|_S$.
 Consequently $[F\,|\,t\cdot I_k]$ is contained in $B$. This is a contradiction so we conclude that $B \in \avbh(m',F)$. 
 Hence $\ncols{B}\le\bsf(m',F)\le \bsf(m,F)$ (by \lrf{monotonicity}). We also know that $\ncols{B} \geq \ncols{A}  - (tk+\ell)m$ 
 and so $\ncols{A} \leq \bsf(m,F) + (tk+\ell)m$ for all $A$.  \qed
\vskip 10pt

\begin{remark}\label{rowcol0s}Let $F$ be a given $k$-rowed (0,1)-matrix.  Let $F'$ denote the matrix obtained from $F$ by adding a row of 0's. Then
$\bsf(m,F')=\bsf(m,F)$ for $m>k$. Also $\bsf(m,[\0_k\,F])=\max\{\ncols{F}, \bsf(m,F))\}$.\end{remark}

\proof Let $A$ be a simple $m$-rowed matrix with $\ncols{A}>\bsf(m,F)$. Then   $F\bs A$.  Now as long as $m\ge k+1$ we have that 
$F'\bs A$.  Similarly if $\ncols{A}>\ncols{F}$, then  $[\0_k\,F]\bs A$. \qed

\vskip 10pt
A more general result would be the following.

\begin{thm}Let $F_1,F_2$ be given. For $F$ as below, $\bsf(m,F)$ is \hfil\break
$O(\ncols{F_1}+\ncols{F_2}+\max\{\bsf(m,F_1),\bsf(m,F_2)\})$.
$$F=\left[\begin{array}{@{}c@{}|@{}c@{}}\begin{array}{c}F_1\\ \hline 0\\ \end{array}&\begin{array}{c}0\\ \hline F_2\\ \end{array}\\ \end{array}\right].$$\label{F1F2}\end{thm}
\proof Assume $F_1$ is $k$-rowed.  Let $A\in\avbh(m,F)$. If $\ncols{A}>\bsf(m,F_1)$, then $F_1\bs A$. Assume $F_1$ appears in the first $k$ rows so that 
$$A=\left[\begin{array}{@{}c@{}|@{}c@{}}\begin{array}{c}F_1\\ \hline *\\ \end{array}&\begin{array}{c}*\\ \hline B\\ \end{array}\\ \end{array}\right].$$
If $F_2\bs B$ then $F\bs A$ and so we may assume $F_2\not\llcurly B$. Now the multiplicity of any column of $B$ is at most $2^k$. 
Thus $\ncols{B}\le 2^k\bsf(m,F_2)$ and so $\ncols{A}\le \ncols{F_1}+ 2^k\bsf(m-k,F_2)\le \ncols{F_1}+ 2^k\bsf(m,F_2)$ by \lrf{monotonicity}. Interchanging $F_1,F_2$ yields the result.
\qed

\section{$3\times \ell$  Berge hypergraphs}\label{k=3}
This section provides an explcit classification of the asymptotic bounds $\bsf(m,F)$. Let
$$
G_1=\left[\begin{array}{cc}1&1\\ 1&0\\ 0&1\\ \end{array}\right],\qquad
G_2=\left[\begin{array}{ccc}1&1&0\\ 1&0&1\\ 0&1&1\\ \end{array}\right].$$

\begin{thm}\label{classifyk=3}  Let $F$ be a $3\times \l$ (0,1)-matrix.\\
(Constant Cases) If $F\bs [ I_3\,|\,t\cdot\0_3]$, then $\bsf(m,F)$ is $\Theta(1)$.\\
(Linear Cases) If $F$ has a Berge hypergraph $2\cdot\1_1$ or $\1_2$ and if $F\bs [G_1\, t\cdot [\0\,|\,I_3]]=H(1,3,t)$  
 then $\bsf(m,F)=\Theta(m)$.\\
(Quadratic Cases) If $F$ has a Berge hypergraph $2\cdot \1_2$ or  $G_2$, or $\1_3$ and if  
$F\bs [\1_3\,|\,t\cdot G_2]=H(2,3,t)$ for some $t$, then  $\bsf(m,F)=\Theta(m^2)$.\\
(Cubic Cases) If $F$  has a Berge hypergraph $2\cdot\1_3$  then  $\bsf(m,F)=\Theta(m^3)$.\\
\end{thm}

\proof  The lower bounds follow from \lrf{clique} and \lrf{construction}. 

The constant upper bound for $[ I_3\,|\,t\cdot\0_3]$ is given by \trf{Ik} combined with \lrf{rowcol0s} to add columns of 0's. 
An exact linear bound for $G_1$ is in \trf{G1}. The linear bound for $[G_1\,t\cdot [\0\,|\,I_3]]=H(1,3,t)$  and the 
 quadratic upper bound for  $[\1_3\,|\,t\cdot G_2]=H(2,3,t)$  follow from \trf{1pxIk-p}.  The cubic upper bound for $t\cdot K_3$ follows from \trf{block1s}. 

To verify that all $3$-rowed matrices are handled we first note that $\bsf(m,2\cdot\1_3)$ is $\Theta(m^3)$.  Consider matrices $F$ with $2\cdot\1_3\not\bs F$. Then $F\bs H(2,3,t)$ and so
$\bsf(m,F)$ is $O(m^2)$.  If $2\cdot \1_2$, $\1_3$ or $G_2 \bs F$ then $\bsf(m,F)$ is $\Omega(m^2)$. 
Now assume  $2\cdot \1_2$, $\1_3$ or $G_2 \not\bs F$. Then $G(F)$ (from \drf{G(F)})  has no 3-cycle nor a repeated edge and so $F\bs H(1,3,t)$. Then $\bsf(m,F)$ is $O(m)$.  If $2\cdot \1_1$ or  
$\1_2\bs F$ then $\bsf(m,F)$ is $\Omega(m)$. The only 3-rowed $F$ with $2\cdot \1_1\not\bs F$ and 
$\1_2\not\bs F$ satisfies $F\bs[ I_3\,|\,t\cdot\0_3]$. 
 \qed
\vskip 10pt
The following theorem is an example of the difference between Berge hypergraphs and configurations. Note that $\forb(m,G_1)=2m$ \cite{survey}.

\begin{thm}
	\label{G1}
	$\bsf(m,G_1) = \lfloor \tfrac{3}{2}m \rfloor + 1$
\end{thm}
\proof
 Let $A\in \bAv(m,F)$. Then $A$ has at most $m+1$ columns of sum 0 or 1.  Consider two columns of $A$ of column sum at least 2. If there is a row that has 1's in both column $i$ and column $j$ then we find a Berge hypergraph $G_1$.  Thus columns of column sum at least 2 must occupy disjoint sets of rows and so there are at most $\lfloor\frac{m}{2}\rfloor$ columns of column sum at least 2. This yields the bound. Then we can form an  $A\in \bAv(m,F)$ with $\ncols{A}=\lfloor \tfrac{3}{2}m \rfloor + 1$.  \qed

\section{$4\times \ell$  Berge hypergraphs}\label{k=4}
  Given a (0,1)-matrix  $F$, we denote by  $r(F)$ (the reduction of $F$) the submatrix obtained by deleting all columns of column sum 0 or 1.   In view of \trf{reduction}, we have that $\bsf(m,F)$ is $O(\bsf(m,r(F)))$. 
On 4 rows, there is an interesting and perhaps unexpected result.

\begin{thm}\cite{AKRS}\label{twoconfigs} $\forb(m,\left\{I_2\times I_2,T_2\times T_2\right\})$ is $\Theta(m^{3/2})$ where
$$I_2\times I_2=C_4=\left[\begin{array}{cccc}
1&1&0&0\\ 0&0&1&1\\ 1&0&1&0\\ 0&1&0&1\\\end{array}\right],\,\,
T_2\times T_2=\left[\begin{array}{cccc}
1&1&1&1\\ 0&0&1&1\\ 1&1&1&1\\ 0&1&0&1\\\end{array}\right].$$
\end{thm}

The above result uses the lower bound construction (projective planes).
from the much cited paper of K\H ovari, S\'os and Tur\'an.   
\begin{thm}\cite{KST} $f(C_4,I_{m/2}\times I_{m/2})$ is $\Theta(m^{3/2})$. \label{KSTlower}\end{thm}

We conclude a Berge hypergraph result much in the spirit of Gerbner and Palmer \cite{GP}.  They maximized a different  
extremal function: essentially the number of 1's in a matrix in $\avbh(m,C_4)$. 

\begin{thm}$\bsf(m,C_4)$ is $\Theta(m^{3/2})$\label{C4}\end{thm}
\proof The lower bound follows from \cite{KST}.  It is straightforward to see that $C_4\bs T_2\times T_2$ and then we apply \trf{twoconfigs} for the upper bound. \qed

\vskip 10pt
We give an alternative  argument in \srf{graph} that handles $F=I_2\times I_s$ for $s\ge 2$.  Other 4-rowed Berge hypergraph cases are more straightforward. Let
$$
H_1=\left[\begin{array}{ccc}1&0&0\\ 1&1&0\\ 0&1&1\\ 0&0&1\\  \end{array}\right],\,\,
H_2=\left[\begin{array}{ccc}1&1&1\\ 1&0&0\\ 0&1&0\\ 0&0&1\\ \end{array}\right],\,\,
H_3=\left[\begin{array}{cc}1&1\\ 1&1\\ 1&0\\ 0&1\\  \end{array}\right],\,\,
H_4=\left[\begin{array}{ccccc}1&1&1&0&0\\ 1&0&0&1&1\\ 0&1&0&1&0\\ 0&0&1&0&1\\  \end{array}\right],$$
$$H_5 = \left[\begin{array}{ccc} 1 & 1 & 0 \\ 1 & 1 & 1 \\ 1 & 0 & 1 \\ 0 & 1 & 1  \end{array}\right],\quad
	H_6 = \left[\begin{array}{ccc} 1 & 1 & 0 \\ 1 & 1 & 0 \\ 1 & 0 & 1 \\ 0 & 1 & 1 \end{array}\right],\quad 
	H_7 = \left[\begin{array}{cccc}  1 & 1 & 0 & 0 \\ 1 & 0 & 1 & 0 \\ 1 & 0 & 0 & 1 \\ 0 & 1 & 1 & 1 \end{array}\right].$$

\begin{thm}\label{classifyk=4}  Let $F$ be a $4\times \l$ (0,1)-matrix.\\
(Constant Cases) If $F\bs [I_4\,|\, t\cdot\0_4]$, then $\bsf(m,F)$ is $\Theta(1)$.\\
(Linear Cases) If $F$ has a Berge hypergraph $2\cdot \1_1$ or $\1_2$ and if $r(F)$ is a configuration in $H_1$ or $H_2$ 
 then $\bsf(m,F)=\Theta(m)$.\\
(Subquadratic Cases) If $r(F)$ is $C_4$, then $\bsf(m,F)$ is $\Theta(m^{3/2})$.\\
(Quadratic Cases) If $F$ has a Berge hypergraph $2\cdot \1_2$ or $G_2$, or $\1_3$ and if  $F\bs H(2,4,t)$  for some $t$, then  $\bsf(m,F)=\Theta(m^2)$.\\
(Cubic Cases) If $F$  has a Berge hypergraph $2\cdot\1_3$ or $\1_4$ or  $K_4^2$  or $H_6$ or $H_7$ and if 
$F\bs H(3,4,t)$ then  $\bsf(m,F)=\Theta(m^3)$.\\
(Quartic Cases)  If $F$  has a Berge hypergraph $2\cdot\1_4$  then  $\bsf(m,F)=\Theta(m^4)$.\\
\end{thm}

\proof  The lower bounds follow from \lrf{clique}, \lrf{construction} and \trf{KSTlower}.

The constant upper bound for $[ I_4\,|\,t\cdot\0_4]$ is given by \trf{Ik}  combined with \lrf{rowcol0s} to add columns of 0's. 
The linear upper bound for $F$ where $G(F)$ is a tree (or forest) follows from \trf{tree}. There are only two trees on 4 vertices namely $H_1$ and $H_2$.  Note $[H_2\,|\,t\cdot [\0_4\,|\,I_4]]=H(1,4,t)$. Thus $\bsf(m, [H_2\,|\,t\cdot [\0_4\,|\,I_4]] )$ is $O(m)$ by  \trf{1pxIk-p}. 
Also $\bsf(m, [H_1\,|\,t\cdot [\0_4\,|\,I_4]] )$ is $O(m)$ by Reduction \lrf{morereduction}.  Now \trf{C4} establishes $\bsf(m,C_4)$. 
 The quadratic upper bound for  $H(2,4,t)$   and the   cubic upper bound for $H(3,4,t)$ follow from \trf{1pxIk-p}.  The quartic upper bound for $t\cdot K_4$ follows from \trf{block1s}. 

To verify that all $4$-rowed matrices are handled we first note that $\bsf(m,2\cdot\1_4)$ is $\Theta(m^4)$.  Consider matrices $F$ with $2\cdot\1_4\not\bs F$. Then $F\bs H(3,4,t)$ and so
$\bsf(m,F)$ is $O(m^3)$.   If $2\cdot \1_3\bs F$, then $\bsf(m,F)$ is $\Omega(m^3)$ by \lrf{construction}.  If   $\1_4$,   $K_4^2$,  $H_6$ or $H_7$  $\bs F$ then $\omega(G(F))=4$  and so $\bsf(m,F)$ is $\Omega(m^3)$ by \lrf{clique}. 
  
   The column minimal simple (0,1)-matrices $F$  with $\omega(G(F))=4$ and with column sums at least 2 are   $\1_4$,   $K_4^2$,  $H_5$, $H_6$ and $H_7$. Since $H_6\bs H_5$ it suffices to drop $H_5$ from the list. Now assume $\omega(G(F))\le 3$ and so
$\1_4$,   $K_4^2$,  $H_6$ or $H_7\not\bs F$. Also assume $2\cdot \1_3\not\bs F$.  Let $3,4$ be the rows so that no column has 1's in both rows $3,4$.  Three columns of sum 3 in $F$  either force $\omega(G(F))=4$ or we have a column of sum 3 repeated. So $F$ has at most 2 (different) columns of sum 3 and  $F\bs H(2,4,t)$.

Now assume  $F\bs H(2,4,t)$ but $2\cdot \1_2$, $\1_3$ or $G_2 \not\bs F$. Then $G(F)$ (from \drf{G(F)})  has no 3-cycle nor a repeated edge and so $G(F)$ is a subgraph of $K_{2,2}$ or $K_{1,3}$. In the latter case, $F\bs H(1,4,t)$. Then $\bsf(m,F)$ is $O(m)$.  In the former case, $F\bs H((2,2),t)$ and so \trf{C4} applies.

If $2\cdot \1_1$ or  
$\1_2\bs F$ then $\bsf(m,F)$ is $\Omega(m)$. The only 4-rowed $F$ with $2\cdot \1_1\not\bs F$ and 
$\1_2\not\bs F$ satisfies $F\bs [I_4\,|\,t\cdot\0_4]$. 
 \qed
\vskip 10pt

We give  some exact linear bounds.
$$\hbox{ Let }H_8=\left[\begin{array}{cc}1&0\\ 1&0\\ 0&1\\ 0&1\\ \end{array}\right]$$
 For the following you may note that $\forb(m,F)$ is $\binom{m}{2}+2m-1$ \cite{survey}.

\begin{thm}
	$\bsf(m,H_8) = 2m$.	
\end{thm}

\proof
Let $A \in \bAv(m,H_8)$. Assume that $A$ is a downset by \lrf{shifting}.  Let $A'=r(A)$.  Since $H_8$ has column sums 2 then 
$\bsf(m,H_8)\le\ncols{A'}+m+1$.   If $A'$ has a column of column sum 4 (or more), then $H_8\bs A'$ since $H_8$ has only 4 rows and is simple. If $A'$ has a column of sum 3 say with 1's on rows 1,2,3, then we find $[K_3^3K_3^2]$
in those 3 rows.   If $A'$ has a column of column sum 3, say with 1's in rows 1,2,3 then we cannot have 
a column with a 1 in row 1 and a 1 in row 4 else $H_8\bs A'$ using the fact that $A$ is a downset (using the columns with 1's in rows 1,4 and the column with 1's in rows 2,3). If $A'$ has only columns of  sum 2 then    we deduce that $\ncols{A'}\le m-1$ and so $\bsf(m,H_8)\le 2m$.

The construction to achieve the bound is to take the $m-1$  columns of sum 2 that have a 1 in row 1 as well as all columns of sum 0 or 1. We conclude 
that  $\bsf(m,H_8)=2m$.
\qed

\begin{thm}
	$\bsf(m,H_2) \leq 4\lfloor m/3\rfloor+m+1$.
\end{thm}

\proof Proceed as above. Let $A \in \bAv(m,H_2)$. Assume that $A$ is a downset by \lrf{shifting}.  Let $A'$=$r(A)$ then $\bsf(m,H_2)\le\ncols{A'}+m+1$ since $H_2$ has column sum 2.  If $A'$ has a column of column sum 4 (or more), then $H_2\bs A'$ since $H_2$ has only 4 rows and is simple.  If $A'$ has a column of sum 3 say with 1's on rows 1,2,3, then we find 
$[K_3^3K_3^2]$ 
in those 3 rows.  If $A'$ has such a column of column sum 3,  then $A'$ cannot have 
a column with a 1 in row 1 and a 1 in row 4 else $F\bs A$ using the fact that $A$ is a downset (using the columns with 1's in rows 1,2 and the column with 1's in rows 1,3 and the column with 1's in rows 1,4).  Thus the number of columns of sum 3 is at most $\lfloor m/3\rfloor$.

Let $t$ be the number of columns of sum 3.  If $m=3t$, then we can include all columns of sum 2 that are in the downset of the columns of sum 3. All other columns of sum 2 have their 1's in the $m-3t$ rows disjoint from those of the 1's in the columns of sum 3.  The columns of sum 2, when interpreted as a graph, cannot have a vertex  of degree 3 else $H_2\bs A$.  So the number of columns of sum 2 is at most 
$m-3t$ for $m-3t\ge 3$ and 0 otherwise.  This yields an upper bound.

 A construction to achieve our bound is to simply take  $\lfloor m/3\rfloor$ columns of sum 3 each having their 1's on disjoint sets of rows and then, for each column of sum 3, add 3 columns of sum 2 whose 1's lie in the rows occupied by the 1's of the column of sum 3. 
\

\section{$5\times \ell$  Berge hypergraphs}\label{k=5}

First we give the 5-rowed classification which requires \corf{conjC4} to be true.
\begin{thm}\label{classifyk=5}  Let $F$ be a $5\times \l$ (0,1)-matrix.
Assume $\bsf(m,\1_1\times C_4)$ is $\Theta(m^2)$.\\
(Constant Cases) If $F\bs [I_5\,|\, t\cdot\0_5]$, then $\bsf(m,F)$ is $\Theta(1)$.\\
(Linear Cases) If $F$ has a Berge hypergraph $\1_2$ or $[1\,1]$  and if $r(F)$ is a vertex-edge incidence matrix of a tree 
 then $\bsf(m,F)=\Theta(m)$.\\
(Subquadratic Cases) If $r(F)$ is  is a vertex-edge incidence matrix  of a bipartite graph $G$ with a cycle  then $\bsf(m,F)$ is $\Theta(\ex(m,G))$ i.e. $\Theta(m^{3/2})$.\\
(Quadratic Cases) If $F$ has a Berge hypergraph $2\cdot \1_2$  or $\chi(G(F))\ge 3$, and if  $r(F)$ is a configuration in
$H(2,5,t)$ from \rf{H(p,k,t)} for some $t$ or in $H((1,2,2),t)$ from \rf{generalH} , then  $\bsf(m,F)=\Theta(m^2)$.\\
(Cubic Cases) If $F$  has a Berge hypergraph $2\cdot\1_3$ or $\1_4$  or $K_4^2$  or $H_6$ or $H_7$ and if 
$F\bs H(3,5,t)$ from \rf{H(p,k,t)} for some $t$ then  $\bsf(m,F)=\Theta(m^3)$.\\
(Quartic Cases) If $F$  has a Berge hypergraph $2\cdot\1_4$   or if $\omega(G(F))=5$  and $F\bs H(4,5,t)$ then  $\bsf(m,F)=\Theta(m^4)$.\\
(Quintic Cases)  If $F$  has a Berge hypergraph $2\cdot\1_5$  then  $\bsf(m,F)=\Theta(m^5)$.\\
\end{thm}

\proof  The lower bounds follow from \lrf{clique}, \lrf{construction} and also  \trf{KSTlower}. Note that a bipartite graph on 5 vertices with a cycle must have a 4-cycle.  In the quadratic cases, we could have listed three minimal examples of Berge hypergraphs with
$\chi(G(F))\ge 3$, namely $\1_3$, $G_2$ or the $5\times 5$ vertex edge incidence matrix of the 5-cycle.

The constant upper bound for $[ I_5\,|\,t\cdot\0_5]$ is given by \trf{Ik}  combined with \lrf{rowcol0s} to add columns of 0's. 
The linear upper bound for $F$ where $G(F)$ is a tree (or forest) follows from \trf{tree}. There are a number of trees on 5 vertices.  Let $F$ be the vertex-edge incidence matrix of a bipartite graph  on 5 vertices that contains a cycle and hence contains $C_4$.  Thus $F\bs I_2\times I_3$ and so \trf{K2t} establishes that
$\bsf(m,F)$ is $O(m^{3/2})$. 
 The quadratic upper bound for  $H(2,5,t)$   and the   cubic upper bound for $H(3,5,t)$ and the quartic upper bound for $H(4,5,t)$ follow from \trf{1pxIk-p}.  The quadratic bound for $H((1,2,2),t)$ is \trf{1xC4} under the assumption 
 $\bsf(m,\1_1\times C_4)$ is $\Theta(m^2)$.
The quintic upper bound for $t\cdot K_5$ follows from \trf{block1s}. 

To verify that all $5$-rowed matrices are handled we first note that $\bsf(m,2\cdot\1_5)$ is $\Theta(m^5)$.  Consider matrices $F$ with $2\cdot\1_5\not\bs F$. Then $F\bs H(4,5,t)$ and so
$\bsf(m,F)$ is $O(m^4)$.  

 If $2\cdot \1_4\bs F$, then $\bsf(m,F)$ is $\Omega(m^4)$ by \lrf{construction}.  If   $\omega(G(F))=5$ then  $\bsf(m,F)$ is $\Omega(m^4)$ by \lrf{clique}. 
  
    Now assume $\omega(G(F))\le 4$ and $2\cdot \1_4\not\bs F$.  Let $4,5$ be the rows so that no column has 1's in both rows $4,5$.  Three columns of sum 4 in $F$  either force $\omega(G(F)=5$ or we have a column of sum 4 repeated. So $F$ has at most 2 (different) columns of sum 4 and  so $F\bs H(3,5,t)$ for some $t$ which yields that $\bsf(m,F)$ is $O(m^3)$. 
    
 If $2\cdot \1_3\bs F$, then $\bsf(m,F)$ is $\Omega(m^3)$ by \lrf{construction}.  If 
 $\1_4$  or $K_4^2$  or $H_6$ or $H_7\bs F$, then $\omega(G(F))\ge 4$ and  then  
 $\bsf(m,F)$ is $\Omega(m^3)$ by \lrf{clique}.  
 
 Now assume $\omega(G(F))\le 3$ and $2\cdot \1_3\not\bs F$. 
  If $\alpha(G(F))\ge 3$, then by taking rows 3,4,5 to be the rows of an independent set of size 3, we have
  $F\bs H(2,5,t)$ and so $\bsf(m,F)$ is $O(m^2)$.    The maximal  graph on 5 vertices with $\omega(G(F))\le 3$ and 
  $\alpha(G(F))\le 2$ is in fact $G(\1_1\times C_4)$. Thus $F\bs H((1,2,2),t)$ for some $t$ and by assumption $\bsf(m,F)$ is $O(m^2)$.  
  
  Now if $2\cdot \1_2\bs F$, then $\bsf(m,F)$ is $\Omega(m^2)$ by \lrf{construction}.  If 
 $\1_3$  or $K_3^2\bs F$  then $\omega(G(F))\ge 3$ and  then  
 $\bsf(m,F)$ is $\Omega(m^2)$ by \lrf{clique}.  Now assume $\omega(G(F))\le 2$ and $2\cdot \1_2\not\bs F$. 
 Thus the columns of $F$ of sum at least 2 must have column sum 2 and there are no repeats of columns of sum 2.  
 The graph $G(F)$ has no triangle. If it is not bipartite then $\chi(G(F))\ge 3$ and then $\bsf(m,F)$ is $\Omega(m^2)$. 
 
 Now assume $2\cdot \1_2\not\bs F$ and $\chi(G(F))\le 2$ and so the columns of sum 2 of $F$ form a bipartite graph $G(F)$ and there are no columns of larger sum.
 The graph $G(F)$ is either a tree in 
 which case $\bsf(m,F)$ is $O(m)$ by \trf{tree} or if there is a cycle it must be $C_4$ and so $\bsf(m,F)$ is $\Omega(m^{3/2})$. But 
 $G(F)$ is a subgraph of $K_{2,3}$ and so we may apply \trf{K2t} (and \trf{reduction}) to obtain $\bsf(m,F)$ is $\Omega(m^{3/2})$.

If $2\cdot \1$ or  
$\1_2\bs F$ then $\bsf(m,F)$ is $\Omega(m)$. The only $F$ with $2\cdot \1_1\not\bs F$ and 
$\1_2\not\bs F$ satisfies $F\bs [I_5\,|\,t\cdot\0_5]$ for some $t$. 
 \qed

If we attempted the classification for 6-rowed $F$ then we would need bounds such as $\bsf(m,I_1\times I_2\times I_3)$ and 
$\bsf(m,I_2\times I_2\times I_2)$. 

\section{Berge hypergraphs from graphs}\label{graph}

 Let $G$ be a graph and let $F$ be the vertex-edge incidence graph. This sections explores some connections for Berge hypergraphs $F$ with  extremal graph theory results.  The first results provides a strong  connection with $\ex(m,K_{s,t})$ and the related problem
 $\ex(m,T,H)$ (the maximum number of subgraphs $T$ in an $H$-free graph on $m$ vertices).  Then we consider the case $G$ is a tree (or forest).  Finally we connect the largest clique size $\omega(G)$ with $\bsf(m,F)$.

 \begin{thm}\label{K2t} Let $F=I_2\times I_t$ be the vertex-edge incidence matrix of the complete bipartite graph $K_{2,t}$.  Then
 $\bsf(m,F)$ is $\Theta(\ex(m,K_{2,t}))$ which is $\Theta(m^{3/2})$. \end{thm}
 \proof It is immediate that  $\bsf(m,F)$ is $\Omega(\ex(m,K_{2,t}))$ since the vertex-edge incidence matrix $A$ of a graph on $m$ vertices  with no subgraph $K_{2,t}$ has $A\in\avbh(m,F)$.  
 
 Now  consider 
 $A\in\avbh(m,F)$. 
 Applying \lrf{shifting}, assume $T_i(A)=A$ for all $i$ and so, when columns are viewed as sets, the columns form a downset.
 Thus for every column $\gamma$ of $A$ of column sum $r$, we have that there are all $2^r$ columns $\alpha$ in $A$ with $\alpha\le\gamma$.  Assume for some column $\alpha$ of $A$  of sum 2  that there are  $2^{t-1}$ columns $\gamma$ of $A$ with
 $\alpha\le \gamma$.  But the resulting set of columns have the Berge hypergraph   $\1_2\times I_t$ by \trf{Ik} and then, using the downset idea, will contain the Berge hypergraph $F$. Thus 
  for a given column $\alpha$ of sum 2, there will be at most $2^{t-1}-1$ columns $\gamma$ of $A$ with
 $\alpha< \gamma$. Thus $\ncols{A}\le (2^{t-1})p$ where $p$ is the number of columns of sum 2 in $A$.  We have $p\le\ex(m,K_{2,t})$ which proves the upper bound for $\bsf(m,F)$. 
 \qed
 \vskip 10pt
 Results of Alon and Shikhelman \cite{AlSh} are surprisingly  helpful here. They prove very accurate bounds.
 For fixed graphs $T$ and $H$, let $\ex(m,T,H)$ denote the maximum number of subgraphs $T$ in an $H$-free graph on $m$ vertices.
 Thus $\ex(m,K_2,H)=\ex(m,H)$. The following is  their Lemma 4.4.  The lower bound for $s=3$ can actually be obtained from the construction of Brown \cite{B}. The lower bounds for larger $s$ have also been obtained by Kostochka, Mubayi and Verstra\"ette \cite{KMV}. 
 \begin{lemma}\label{K3inKst}\cite{AlSh} For any fixed $s\ge 2$ and $t\ge (s-1)!+1$, $\ex(m,K_3,K_{s,t})$ is $\Theta(m^{3-(3/s)})$.  \end{lemma}
 We can use this directly in analogy to \trf{K2t}.
 \begin{thm}$\bsf(m,I_3\times I_t)$ is $\Theta(m^{2})$. \label{I3Ik}\end{thm}
 \proof Let 
 $A\in\avbh(m,I_3\times I_t)$. 
 Applying \lrf{shifting}, assume $T_i(A)=A$ for all $i$ and so, when columns are viewed as sets, the columns form a downset. 
 Thus for every column $\gamma$ of $A$ of column sum $r$, we have that there are all $2^r$ columns $\alpha$ in $A$ with $\alpha\le\gamma$.  Let $G$ be the graph associated with the columns of sum 2 and so a column of sum $r$ corresponds to $K_r$ in $G$. In particular the number of columns of sum 3 is bounded by $\ex(m,K_3,K_{3,t})$ since each column of sum 3 yields a triangle $K_3$. 
 Assume for some column $\alpha$ of $A$  of sum 3  that there are  $2^{t-1}$ columns $\gamma$ of $A$ with
 $\alpha\le \gamma$.  But the resulting set of columns have the Berge hypergraph   $\1_3\times I_t$ by \trf{Ik} and then, using the downset idea, will contain the Berge hypergraph $I_3\times I_t$. Thus 
  for a given column $\alpha$ of sum 3, there will be at most $2^{t-1}-1$ columns $\gamma$ of $A$ with
 $\alpha< \gamma$. Thus $\ncols{A}\le (2^{t-1})p+|E(G)|$ where $p$ is the number of columns of sum 3 in $A$.  We have $p\le\ex(m,K_3,K_{3,t})$. 
 This yields $\ncols{A}\le 2^{t-1} \ex(m,K_3,K_{3,t})+ \ex(m,K_{3,t})$. Now the standard inequalities yield $\ex(m,K_{3,t})$ is $O(m^{5/3})$ and combined with \lrf{K3inKst} we obtain the upper bound.  The lower bound would follow from taking construction of 
 $\Theta(m^{3-(3/t)})$ triples as columns of sum 3 from \lrf{K3inKst}.
  \qed
  \vskip 10pt
    We could follow the above proof technique and verify, for example, that 
  $$\bsf(m,I_4\times I_4)\hbox{ is }O\left( \ex(m,K_{4,4})+\ex(m,K_3,K_{4,4})+\ex(m,K_4,K_{4,4})\right)$$
  using the idea that we can restrict our attention, for an asymptotic bound, to columns of sum 2,3,4.
  Note that \lrf{K3inKst} yields $\ex(m,K_3,K_{4,4})$ is  $\Theta(m^{2+(1/4)})$ and so 
  $\bsf(m,I_4\times I_4)$ is $\Omega(m^{2+(1/4)})$. Thus
  $I_{m/2}\times I_{m/2}$ won't be the source of the construction.  The paper \cite{AS} has some lower bounds (Lemma 4.3 in \cite{AS}): 
  \begin{lemma} \cite{AS} For any fixed $r$, $s\ge 2r-2$ and $t\ge (s-1)!+1$. Then 
  $$\ex(m,K_r,K_{s,t})\ge \left(\frac{1}{r!}+o(1)\right)m^{r-\frac{r(r-1)}{2s}}.\qed$$\end{lemma} 
  Thus for some choices $r,s,t$,  $\ex(m,K_r,K_{s,t})$ grows  something like $\Omega(m^{r-\epsilon})$  which shows we can take many columns of sum $r$ and still avoid $K_{s,t}$, i.e.  
  $\bsf(m,K_{s,t})$  grows very large.

 \begin{thm} Let $F$ be the vertex-edge incidence $k\times (k-1)$ matrix of a tree (or forest) $T$ on $k$ vertices.  
 Then $\bsf(m,F)$ is $\Theta(m)$.\label{tree}\end{thm}
\proof We generalize the result for trees/forests in graphs. It is known that if a graph $G$ has all vertices of degree $k-1$, then $G$ contains any tree/forest on $k$ vertices as a subgraph.  We follow that argument but need to adapt the ideas to Berge hypergraphs.  
 Let $A\in\avbh(m,F)$  with $A$ being a downset. We will show that $\ncols{A}\le 2^{k-1}m$.
If $A$ has all rows sums at least $2^{k-1}+1$ then we can establish the result as follows.  If we consider the submatrix $A_r$ formed by those columns with a 1 in row $r$, then $I_{k-1}$ is a Berge hypergraph contained in the rows $[m]\backslash r$ of $A_r$ (by \trf{Ik}). Thus the vertex corresponding to row $r$ in $G(A)$ has degree at least $k-1$.     Then $G(A)$ has a copy of the tree/forest $T$ and since $A$ is a downset, $F\bs A$, a contradiction.

If $A$ has some rows of sum at most $2^{k-1}$, then we use 
 induction on $m$.   
 Assume   row $r$ of $A$  has  row sum 
$t\le 2^{k-1}$. Then we may delete that row and the $t$ columns with 1's in row $r$ and the resulting $(m-1)$-rowed matrix $A'$ is simple with $\ncols{A}=\ncols{A'}+t$. By induction $\ncols{A'}\le 2^{k-1}(m-1)$ and this yields $\ncols{A}\le  2^{k-1}m$.
\qed
 \vskip 10pt
 
 This argument does not extend to other graphs since when we find our desired Berge hypergraph $I_{\ell}$ under 1's on row $r$, we cannot control which rows of $A$ contain a Berge hypergraph $I_{\ell}$. The following results shows the large gap between Berge hypergraph results and forbidden configurations results.
 
 The following matrices will be used in our arguments.
  \begin{equation}
  F_7 = \begin{bmatrix} 
    1 & 1 & 0 & 1 & 1 & 0  \\
    1 & 0 & 1 & 1 & 1 & 1  \\
    0 & 1 & 0 & 1 & 0 & 1  \\
    0 & 0 & 1 & 0 & 0 & 1  \\
    0 & 0 & 0 & 0 & 1 & 0  \\
  \end{bmatrix},\quad 
  H_9=\begin{bmatrix} 1 & 0 & 0 \\ 1 & 0 & 0 \\ 0 & 1 & 0 \\0 & 1 & 0 \\ 0 & 0 & 1 \\  0 & 0 & 1 \\ \end{bmatrix}, quad
  H_{10}=\begin{bmatrix} 1 & 1 & 0  \\  1 & 0 & 0 \\ 0 & 1 & 0 \\ 0 & 0 & 1 \\ 0 & 0 & 1 \end{bmatrix}.\label{F7matrix}\end{equation}
  
 \begin{lemma}
\vspace{1pt}
For $k\geq 5$, $\forb(m,H_1\times \0_{k-4})$ is $\Theta(m^{k-3})$.
\end{lemma}

\proof
The survey \cite{survey} has the result  $\forb(m,F_7)$ is $\Theta(m^2)$ listed  in the results on 5-rowed $F$. We have $H_1\times \0_1 \prec F_7$. Thus $\forb(m,H_1\times \0_1)$ is $O(m^2)$. The upper bound for $k\geq 6$ follows by standard induction. We note that $H_1\times \0_{k-4}$ has a $(k-2)\times l$ submatrix with $K_2^0$ on every pair of rows and so $\forb(m,H_1\times \0_{k-4})$ is $\Omega(m^{k-3})$ by  \cite{AFl}. \qed

\begin{lemma}
$\forb(m,H_2 \times \0_{k-4})$ is $\Theta(m^{k-2})$.
\end{lemma}

\proof
$\forb(m,H_2)$ is $\Theta(m^2)$ from Theorem 6.1 of \cite{survey} so by induction $\forb(m,H_2\times \0_{k-4})$ is $O(m^{k-2})$. $H_2 \times \0_{k-4}$ has a $(k-1)\times l$ submatrix with $K_2^0$ on every pair of rows so by \cite{AFl}, $H_2 \times \0_{k-4}$ is $\Omega(m^{k-2})$. \qed

\begin{lemma}
$\forb(m,H_9\times \0_{k-6})$ is $\Theta(m^{k-1})$.
\end{lemma}
\proof
$H_9\times \0_{k-6}$ has $K_2^0$ on every pair of rows so by \cite{AFl}, $H_9\times \0_{k-6}$ is $\Theta(m^{k-1})$. \qed

\begin{thm}\label{treeforb}
Assume $k\geq 5$ and let $F$ be the $k\times l$ vertex-edge incidence matrix of a forest $T$. 
\begin{enumerate}[i.]
	\item
$\forb(m,F)$ is $\Theta(m^{k-3})$ if and only if $F\prec H_1$.

\item
$\forb(m,F)$ is $\Theta(m^{k-2})$ if and only if $F\nprec H_1$ and $T$ has at most 2 stars connected by a path of at most 2 or not connected.

\item
If $F$ is not one of the two previous cases, then $\forb(m,F) $ is $\Theta(m^{k-1})$.

\end{enumerate}
\end{thm} 

\proof

\noindent
Assume $k\geq 5$.

\noindent
{\bf Case 1: $\forb(m,F)$ is $\Theta(m^{k-3})$.}

Note that $\forb(m,F)$ is $\Omega(m^{k-3})$ for all trees since a single edge produces a column which has $k-2$ rows with $K_2^0$ on every pair of rows. If $T$ has only 2 edges then $F\prec H_1\times 0_{k-4}$ and $\forb(m,F)$ is $\theta(m^{k-3})$. Otherwise if $F$ has 3 columns and $F\nprec H_1\times 0_{k-4}$ and $F \nprec H_2\times \0_{k-4}$ then $F=H_{10} \times \0_{k-5}$ up to isomorphism. Rows 2 to $k$ have the property that every pair of rows has $K_2^0$ so $\forb(m,F)$ is $\Theta(m^{k-2})$. If $F$ has more than three columns then there is a triple of columns, $G$ such that $G\nprec H_1\times \0_{k-4}$ so $\forb(m,F)$ is $\Omega(m^{k-2})$. So $\forb(m,F)$ is $\theta(m^{k-3})$ if and only if $F\prec H_1\times \0_{k-4}$. 

\medskip

\noindent
{\bf Case 2: $F\not\prec H_1$ and $H_9\not\prec F$.}

If $F_3 \nprec F$ then $T$ has no path with 5 or more edges and $T$ has two or fewer non-trivial components. If $T$ has two components then each component is a star since no component has two disjoint edges. Furthermore, if $T$ has a path of 4 edges then the middle vertex has degree 2 and the second and fourth vertex can high degree. If $T$ does not have a path of 4 edges but has a path with 3 edges then the second and third vertices of the path are centers of connected stars. If $T$ has a path of length at most 2 then $T$ is a star. For each possibility of $T$ with $F\nprec F_3$ and $F\nprec H_1 \times \0_{k-4}$, there are a pair of vertices $a$,$b$ such that if $\{a,c\}$ is an edge then $c=b$, a pair of vertices that are not connected, and a pair of vertices $d$,$e$ such that every edge contains either $d$ or $e$. By \trf{smallboundary}, $\forb(m,F)$ is $O(m^{k-2})$. Since $F\nprec H_1\times \0_{k-4}$, $\forb(m,F)$ is $\Omega(m^{k-2})$. This concludes the case $\forb(m,F)$ is $\Theta(m^{k-2})$.

\medskip

\noindent
{\bf Case 3: $H_9\prec F$}

Because $F$ has column sums 2, then $H_9\prec F$ implies  $H_9\times \0_{k-6} \prec F$ and so $\forb(m,F)$ is $\Omega(m^{k-1})$. 
Because $F$ is simple then 
$\forb(m,F)$ is $O(m^{k-1})$ \cite{survey}. 

This concludes the case $k\geq 5$. \qed

\section{Conjecture and Problems}

We have used the following  conjecture in \trf{classifyk=5}.
\begin{conj}\label{conjC4}$\bsf(m,\1_1\times C_4)$ is $\Theta(m^2)$. \end{conj}

 What are the equivalent difficult cases for larger number of rows. The conjecture would yield
 $\bsf(m,\1_2\times C_4)$ is $\Theta(m^3)$ by \lrf{induction}  but we do not predict $\bsf(m,\1_1\times I_2\times I_3)$.
  For $k=6$, we believe that $F=I_2\times I_2\times I_2$ will be quite challenging given an old result of Erd\H{o}s \cite{E64}.

\begin{thm}\cite{E64} $f(I_2\times I_2\times I_2, I_{m/3}\times I_{m/3}\times I_{m/3})$ is $O(m^{11/4})$ and $\Omega(m^{5/2})$.\label{2x2x2}\end{thm}

We might predict that $\bsf(m,I_2\times I_2\times I_2)=\Theta(f(I_2\times I_2\times I_2, I_{m/3}\times I_{m/3}\times I_{m/3}))$.
  and so $\bsf(m,I_2\times I_2\times I_2)$ is between quadratic and cubic.   Unfortunately we offer no improvement to the bounds of
   Erd\H{o}s.

\end{document}